\documentclass{paper}
\usepackage{amssymb,amsmath}
\usepackage{amsfonts}
\usepackage{latexsym,diagrams}
\input{defs_allg.tex.input}

\newarrowhead{3>}{\Rrightarrow}{\Lleftarrow}{\vboxtoz{\vss\hbox{$\vee$}\kern0pt}}{\vboxtoz{\hbox{$\wedge$}\vss}}
\newarrow{Three}{3}{3}{3}{3}{3>}

\begin{document}

\begin{frontmatter}
\title{Nichols algebras of $U_q(\mathfrak{g})$-modules\thanksref{PhD}}
\thanks[PhD]{This work will be part of the authors PhD thesis written under the supervision of Professor H.-J. Schneider.}
\author{Stefan Ufer\thanksref{BayStaat}}
\address{Mathematisches Institut der Universit\"at M\"unchen, Theresienstr. 39,\\ 80333 M\"unchen, Germany}
\thanks[BayStaat]{Partially supported by Graduiertenf\"orderung des bayerischen Staates and by his parents}
\begin{abstract}
A technique is developed to reduce the investigation of Nichols algebras of integrable $U_q(\mathfrak{g})$-modules to the investigation of Nichols algebras of braided vector spaces with diagonal braiding. The results are applied to obtain information on the Gelfan'd-Kirillov dimension of these Nichols algebras and their defining relations if the braiding is of a special type and $q$ is not a root of unity. \\
For simple $\mathfrak{g}$ the cases of finite Gelfand-Kirillov dimension are determined completely.
\end{abstract}
\end{frontmatter}
\section{Introduction}

Given a braiding $c$ on a vector space $V$ one can define the \emph{Nichols algebra} $\Nichols(V,c)$ of $V$ and $c$ as a certain generalization of the classical symmetric algebra (definiton \ref{defn_Nichols}). The symmetric algebra $S(V)$ is obtained if $c=\tau$ is the usual flip map 
\[\tau: V\otimes V\rightarrow V\otimes V, \: \tau(v\otimes w)= w\otimes v.\]
These Nichols algebras occured first in Nichols paper \cite{Nichols}. The aim there was to classify certain finite-dimensional Hopf algebras. For a survey on Nichols algebras see \cite{AS5}.\par
Apart from special - though very important - classes of braidings not very much is known about these algebras. In particular one is interested in the vector space dimension or in the Gelfand-Kirillov dimension and in defining relations.\par
Much is known in the case of diagonal braidings. The Nichols algebras of these braidings play a central role in the theory of quantum groups \cite{Lusztig} (Lusztigs algebra $\mathbf{f}$) and in the classification of pointed Hopf algebras \cite{AS5}. The key idea of this paper is to use the known results for the diagonal case from \cite{AS4,Rosso_invent} to obtain new results for more general braidings. \par
An important source of non-diagonal braidings are finite-dimensional integrable modules over quantum enveloping algebras $U_q(\mathfrak{g})$ of finite-dimensional semisimple complex Lie algebras $\mathfrak{g}$. In \cite{A1} Andruskiewitsch raises the question on the structure of the induced Nichols algebras. Partial results are known in this case, e.g. Rosso \cite{Rosso_invent} determined the structure for certain simple $U_q(\mathfrak{sl}_n)$-modules (in these cases the braidings are of Hecke type). For some low-dimensional simple $U_q(\mathfrak{sl}_2)$-modules (some not of Hecke type) the author determined the structure in \cite{Ich_PBW}.\par
Inspired by an observation of Rosso in \cite{Rosso_invent} this paper describes a method to reduce the study of Nichols algebras of finite-dimensional integrable $U_q(\mathfrak{g})$-modules $M$ to the study of Nichols algebras of diagonal braidings. This method consists of two steps.
The first step is to define a Hopf algebra extension $\hat{U}$ of the negative part of $U_q(\mathfrak{g})$ such that $M$ is a Yetter-Drinfel'd module over $\hat{U}$ and the given braiding on $M$ is equal to the one induced by the Yetter-Drinfel'd module structure. For simple modules $M$ a similar construction is already mentioned in \cite{Rosso_invent}.
The extension $\hat{U}$ decomposes into a Radford biproduct
\[ \hat{U} = \Nichols(V)\#kG,\]
where $V$ is a completely reducible Yetter-Drinfel'd module over the abelian group $G$.
In the second step the general theorem \ref{thm_nicholsext} says that the braided biproduct of $\Nichols(M)$ and $\Nichols(V)$ is again a Nichols algebra of a completely reducible Yetter-Drinfel'd module $M_{kG}\oplus V$ over $G$:
\[ \Nichols(M)\#\Nichols(V) \isom \Nichols(M_{kG}\oplus V).\]
Here $M_{kG}$ is the linear subspace of $M$ spanned by the vectors of highest weight. It becomes a Yetter-Drinfel'd module over $G$ in a natural way, and the braiding is diagonal.
Actually theorem \ref{thm_nicholsext} is proven for arbitrary Hopf algebras $H$ instead of the group algebra $kG$.

The described method is applied to obtain a criterion for the finiteness of the Gelfand-Kirillov dimension of $\Nichols(M)$. To the module $M$ and its braiding one associates a matrix $(b_{ij})$ of rational numbers that is an extension of the Cartan matrix of $\mathfrak{g}$. Under some technical assumptions on the braiding, the Gelfand-Kirillov dimension of $\Nichols(M)$ is finite if and only if $(b_{ij})$ is a Cartan matrix of finite type. For simple $\mathfrak{g}$ and simple modules $M$ a complete list of all cases with finite Gelfand-Kirillov dimension is given in table \ref{table_extensions}.\par

As a second application of the method a description of the relations of $\Nichols(M)$ is obtained under the assumption that the braided biproduct $\Nichols(M)\#\Nichols(V)$ is given by the quantum Serre relations. In particular this applies if $\Nichols(M)$ has finite Gelfand-Kirillov dimension. Table \ref{table_extensions} contains the degrees of the defining relations in the case that $\mathfrak{g}$ is simple, $M$ is simple and the Gelfand-Kirillov dimension of $\Nichols(M)$ is finite.\par
Due to missing information on Nichols algebras of diagonal braidings the results of both applications contain some technical restrictions.\par
The paper is organized as follows. In section \ref{sect_basic} some basic definitions and results are recalled. These include braidings, Yetter-Drinfel'd modules, braided Hopf algebras, Nichols algebras and Radfords theorem for Hopf algebras with a projection. Section \ref{sect_uqgmod} contains the required material on $U_q(\mathfrak{g})$-modules, mainly taken from the book of Jantzen \cite{Jantzen}. In section \ref{sect_ydstruct} the Hopf algebra extension $\hat{U}$ of $U_q^{\leq 0}(\mathfrak{g})$ and a functor from the category of integrable $U_q(\mathfrak{g})$-modules to the category of Yetter-Drinfel'd modules over $\hat{U}$ are defined such that the functor preserves the braiding. 
Section \ref{sect_brbiprod} contains the definition of a braided version of Radfords biproduct and some necessary results on these objects. The last important notion for the extension theorem in section \ref{sect_extthm}, graded Yetter-Drinfel'd modules, is introduced in section \ref{sect_gryd}. Section \ref{sect_uqgapp} contains the results on the Gelfand-Kirillov dimension of the Nichols algebras and in section \ref{sect_rel} results on the defining relations can be found.\par
Assume that $k$ is a field. Tensor products are always taken over $k$. Throughout the paper $q\in k$ is not a root of unity. Sweedlers notation without the sigma-sign will be used for comultiplications and coactions according to the following convention:
\begin{nota}
Assume that $H$ is a Hopf algebra with bijective antipode, $R$ a Hopf algebra in ${}^H_H\mYD$ such that $R\#H$ has bijective antipode and that $Q$ is a braided Hopf algebra in ${}^{R\#H}_{R\#H}\mYD$. In this paper the following conventions for Sweedler notation are used:
\begin{enumerate}
\item The Sweedler indices for the comultiplication in usual Hopf algebras are lower indices with round brackets: $\Delta_H(h) = h\s1\otimes h\s2$.
\item The Sweedler indices for the comultiplication in braided Hopf algebras in ${}^H_H\mYD$ are upper indices with round brackets: $\Delta_R(r) = r\su1\otimes r\su2$.
\item The Sweedler indices for the comultiplication in braided Hopf algebras in ${}^{R\#H}_{R\#H}\mYD$ are upper indices with square brackets: $\Delta_Q(x) = x\sxu1\otimes x\sxu2$.
\item For $H$ coactions use lower Sweedler indices with round brackets: $\delta_H(v) = v\sm1\otimes v\s0$.
\item  For $R\#H$ coactions use lower Sweedler indices with square brackets: $\delta_{R\#H}(m) = m\smx1\otimes m\sx0$.
\end{enumerate}
\end{nota}
\newpage

\section{Basic concepts}
\label{sect_basic}
In this section some basic definitions and facts are listed. This will include braidings, Yetter-Drinfel'd modules, braided Hopf algebras, Radford biproducts and Radfords theorem for Hopf algebras with a projection.
\subsection{Braidings}
\begin{defn}
A \emph{braided vector space $(V,c)$} is a vector space $V$ together with a $k$-linear automorphism $c$ of $V\otimes V$ such that the following identity (the braid equation) holds in $\End_k(V\otimes V\otimes V)$:
\[ (c\otimes\id_V)(\id_V\otimes c)(c\otimes\id_V) = (\id_V\otimes c)(c\otimes\id_V)(\id_V\otimes c).\]
Let $X$ be a basis of $V$.
The braiding is \emph{diagonal with respect to $X$} if a matrix $(q_{xy})_{x,y\in X}$ of non-zero scalars such that
\[ \forall x,y\in X: c(x\otimes y) = q_{xy} y\otimes x.\]
The braiding is called \emph{diagonal} if it is diagonal with respect to some basis of $V$.
\end{defn}
Recall that a \emph{generalized Cartan matrix} \cite{Kac} is a matrix $(a_{xy})_{x,y\in X}$ (where $X$ is a finite index set) with integer coefficients such that
\begin{itemize}
\item $\forall x\in X: a_{xx} = 2$,
\item $\forall x,y\in X,x\neq y: a_{xy} \leq 0$ and
\item $\forall x,y\in X: a_{xy}=0 \Rightarrow a_{yx}=0$.
\end{itemize}
The finite-dimensional braided vector space $(V,c)$ is called \emph{of Cartan type} (with respect to the basis $X$) if it is diagonal (with respect to $X$) and the coefficients $(q_{xy})$ satisfy
\[\forall x,y\in X: q_{xy}q_{yx} = q_{xx}^{a_{xy}}.\]
for a generalized Cartan matrix $(a_{xy})_{x,y\in X}$.\par
To a generalized Cartan matrix $(a_{xy})_{x,y\in X}$ define its Coxeter graph as the (unoriented) graph with vertex set $X$ and $a_{xy}a_{yx}$ edges between $x$ and $y$ for all $x\neq y\in X$. The generalized Cartan matrix is called \emph{connected} if its Coxeter graph is connected. Define the Dynkin diagram of $(a_{xy})$ as in \cite[\S 4.7.]{Kac}. A generalized Cartan matrix is called \emph{of finite type} if it is the Cartan matrix of a finite-dimensional semisimple complex Lie algebra.

\subsection{Yetter-Drinfel'd modules}
Let $H$ be a Hopf algebra with bijective antipode.
\begin{defn}
A (left-left) Yetter-Drinfel'd module $M$ over $H$ is a left $H$-module and a left $H$-comodule such that for all $h\in H,m\in M$
\[ (hm)\sm1\otimes (hm)\s0 =h\s1 m\sm1 S\left(h\s3\right)\otimes h\s2m\s0.\]
\end{defn}
The category of Yetter-Drinfel'd modules ${}^H_H\mYD$ over $H$ is a braided monoidal category with braidings given by
\[ c_{M,N} :M\otimes N\rightarrow N\otimes M, c_{M,N}(m\otimes n) := m\sm1 n\otimes m\s0.\]
In particular every Yetter-Drinfel'd module $M$ is a braided vector space with braiding $c_{M,M}$.

\subsection{Braided Hopf algebras}

Assume that $(R,m,\eta)$ is a $k$ algebra and $R\in{}^H_H\mYD$ such that $m,\eta$ are $H$-linear and $H$-colinear. Then $R\obar R:=R\otimes R$ is an algebra with unit $1_R\otimes 1_R$ and multiplication given by
\[ m_{R\obar R} := (m\otimes m)(\id_R\otimes c_{R,R}\otimes \id_R).\]
\begin{defn}
Let $H$ be a Hopf algebra with bijective antipode. $(R,m,\eta,\Delta,\eps)$ is called a \emph{braided bialgebra in ${}^H_H\mYD$} if
\begin{itemize}
\item $m,\eta,\Delta,\eps$ are $H$-linear and $H$-colinear,
\item $(R,m,\eta)$ is an algebra,
\item $(R,\Delta,\eps)$ is a coalgebra and
\item $\Delta:R\rightarrow R\obar R$, $\eps:R\rightarrow k$ are algebra morphisms.
\end{itemize}
If the identity has an inverse $S$ in the convolution algebra $\End_k(R)$, then $(R,m,\eta,\Delta,\eps)$ is called a \emph{braided Hopf algebra in ${}^H_H\mYD$}. In this case $S$ is called the antipode of $R$. It is $H$-linear and $H$-colinear.
\end{defn}

\subsection{Nichols algebras}
Important examples of braided Hopf algebras are Nichols algebras which are a generalization of the symmetric algebra, replacing the usual flip map with an arbitrary braiding.
\begin{defn}
\label{defn_Nichols}
Let $H$ be a Hopf algebra with bijective antipode, $M$ an ${}^H_H\mYD$-module and $c = c_{M,M}$. The Nichols algebra $\Nichols(M,c)$ is a braided Hopf algebra in ${}^H_H\mYD$ satisfying
\begin{itemize}
\item $\Nichols(M,c) = \bigoplus_{n\geq 0} R(n)$ is graded as an algebra and as a coalgebra,
\item $c(R(i)\otimes R(j))\subset R(j)\otimes R(i)$,
\item $R(1)\isom M$ as braided vector spaces
\item $\Nichols(M,c)$ is generated by $R(1)$ as an algebra and
\item $R(1) = P(R)$.
\end{itemize}
\end{defn}
The Nichols algebra exists for every Yetter-Drinfel'd module $M$ and it is unique up to isomorphism \cite{AS5,Sbg_borel}. More generally one can define the Nichols algebra for any braided vector space. It is then a braided Hopf algebra in the sense of \cite{Takeuchi_survey} and does only depend on the braiding $c$ and not on the Yetter-Drinfel'd structure. In this article only the Yetter-Drinfel'd case will be considered.

\subsection{The Radford biproduct}
Let $H$ be a Hopf algebra with bijective antipode and $R$ a braided Hopf algebra in ${}^H_H\mYD$. Define the Radford biproduct of $R$ with $H$ as the vector space $R\# H:=R\otimes H$. By \cite{Radford_HAproj} it is a Hopf algebra with unit $1_R\# 1_H$, counit $\eps_R\otimes\eps_H$, multiplication
\[ \forall r,r'\in R,h,h'\in H: (r\#h)(r'\#h') := r(h\s1\cdot r')\#h\s2 h'\]
and comultiplication
\[ \forall r\in R,h\in H: \Delta_{R\#H} (r\#h) := r\su1\# r\su2\sm1 h\s1\otimes r\su2\s0\# h\s2.\]
Note that the map
\[ \pi:=\eps\#H: R\#H\rightarrow H\]
is a Hopf algebra projection and $\pi(1\#h)=1\#h$ for all $h\in H$. If $R$ has bijective antipode so has $R\#H$.

\subsection{Hopf algebras with a projection}
On the other hand, every Hopf algebra with bijective antipode that has a projection onto a Hopf subalgebra is a Radford biproduct. Let $A$, $H$ be Hopf algebras with bijective antipode and $\pi:A\rightarrow H$ and $\iota:H\rightarrow A$ be morphisms of Hopf algebras such that
\[ \pi\iota = \id_H.\]
Consider $H$ as a Hopf subalgebra of $A$ via $\iota$.\par
The algebra of right coinvariants with respect to $\pi$ is
\[ R := A^{\co\pi} = \{a\in A|\:(\id\otimes\pi)\Delta(a) = a\otimes 1\}.\]
This algebra $R$ is a Yetter-Drinfel'd module over $H$, where the action is given by the adjoint action of $H$ in $A$. The coaction is given by
\[ \delta(r) := (\pi\otimes\id)\Delta(r).\]
$R$ is a subalgebra of $A$, but in general not a subcoalgebra. Nevertheless with the counit $\eps_A|R$ and comultiplication
\[ \Delta_R(r) := r\s1 \iota S_H\pi(r\s2)\otimes r\s3\]
$R$ becomes a braided Hopf algebra in ${}^H_H\mYD$.\par
Moreover the multiplication of $A$ induces an isomorphism of Hopf algebras
\[ R\# H\rightarrow A, \:r\#h\mapsto r\iota(h).\]
Consider the map
\[\theta_R:A\rightarrow R,\:\theta_R(a) := a\s1\iota S_H\pi(a\s2).\]
Then the inverse of the isomorphism above is given by the map
\[ (\theta_R\#\pi)\Delta:A\rightarrow R\#H.\]
Note however that in general $\theta_R$ is not an algebra morphism.

\section{Braidings on $U_q(\mathfrak{g})$-modules}
\label{sect_uqgmod}
In this section some material from the book of Jantzen \cite{Jantzen} is presented. Assume $\operatorname{char} k=0$. Fix a complex finite-dimensional semisimple Lie algebra $\mathfrak{g}$ with root system $\Phi$ and weight lattice $\Lambda$. Abbreviate $U:=U_q(\mathfrak{g})$. Fix a basis $\Pi$ of the root system and let $(a_{\alpha,\beta})_{\alpha,\beta\in\Pi}$ be the Cartan matrix. The $\Q$ vector space $V := \Q\otimes \Lambda$ has a non-degenerate symmetric bilinear form $(-,-)$ such that for all $\alpha,\beta\in\Pi$
\[ a_{\alpha,\beta} = \frac{2(\alpha,\beta)}{(\alpha,\alpha)}.\]
There is a partial order on $\Z\Phi$ defined by
\[ \forall \mu,\nu\in\Z\Phi: \mu>\nu \Leftrightarrow \mu-\nu\in(\N_0\Pi)\setminus \{0\}.\]
Let $U^+$, resp. $U^-$, be the subalgebra of $U$ generated by the $E_\alpha,\alpha\in\Pi$, resp. by the $F_\alpha,\alpha\in\Pi$. $U$, $U^+$ and $U^-$ admit weight space decompositions. Denote the weight space corresponding to the root $\mu\in\Z\Phi$ by $U_\mu$, resp. $U^+_\mu$, resp. $U^-_\mu$.\par
Following \cite[7.1]{Jantzen} choose for all $\mu\in\Z\Phi,\mu\geq 0$ a basis $u_1^\mu,\ldots,u_{r(\mu)}^\mu$ of $U_\mu^+$ and a dual basis (with respect to the non-degenerate pairing between $U^+$ and $U^-$ defined there) $v_1^\mu,\ldots,v_{r(\mu)}^\mu$ of $U_{-\mu}^-$. Define
\[ \Theta_\mu := \sum\limits_{i=1}^{r(\mu)} v_i^\mu\otimes u_i^\mu\in U_{-\mu}^-\otimes U_\mu^+.\]
In the sequel write formally
\[ \Theta_\mu = \Theta_\mu^-\otimes \Theta_\mu^+, \]
keeping in mind that $\Theta_\mu$ is almost always a sum of tensors.\par

\begin{defn}
Let $\sigma:\Z\Phi\rightarrow \{\pm 1\}$ be a group homomorphism.
A $U$-module $M$ is called integrable of type $\sigma$ if it is the direct sum of its weight spaces:
\[ M = \bigoplus\limits_{\lambda\in\Lambda}\{m\in M| \forall \mu\in\Z\Phi : K_\mu m = \sigma(\mu)q^{(\lambda,\mu)} m\}\]
and the $E_\alpha,F_\alpha$ act locally nilpotently on $M$, i.e. for every $m\in M$ there exists $n\in\N$ such that for all $\alpha\in\Pi: E_\alpha^nm=0=F_\alpha^nm$.
\end{defn}
\begin{rem}
For all group homomorphisms $\sigma:\Z\Phi\rightarrow \{\pm 1\}$ the category of $U$-modules of type $\sigma$ is equivalent to the category of $U$-modules of type $+1$ (the trivial homomorphism) \cite[5.2.]{Jantzen}. This equivalence is in general \emph{not} monoidal. Nevertheless in this paper only $U$-modules of type $+1$ are considered. A $U$-module will be called integrable if it is integrable of type $+1$.
\end{rem}
Choose a function $f:\Lambda\times\Lambda\rightarrow k^\times$ such that for all $\lambda,\mu\in\Lambda,\nu\in\Z\Phi$
\begin{equation}
\label{equ_f}
f(\lambda+\nu,\mu) = q^{-(\nu,\mu)}f(\lambda,\mu)\:\mbox{and}\: f(\lambda,\mu+\nu) = q^{-(\lambda,\nu)}f(\lambda,\mu).
\end{equation}
Now for integrable $U$-modules $M,M'$ define a $U$-linear isomorphism
\[c^f_{M,M'}:M\otimes M'\rightarrow M'\otimes M\]
 as in \cite[7.3]{Jantzen} such that for all $m\in M_\lambda,m'\in M'_{\lambda'}$
\[c^f_{M,M'}(m\otimes m') = f(\lambda',\lambda) \sum\limits_{\mu\geq 0} \Theta_\mu (m'\otimes m).\]
On every triplet of integrable $U$-modules $M,M',M''$ these morphisms satisfy the braid equation
\begin{eqnarray*}
(c^f_{M',M''}\otimes \id_{M})(\id_{M'}&\otimes& c^f_{M,M''})(c^f_{M,M'}\otimes\id_{M''}) = \\ &=& (\id_{M''}\otimes c^f_{M,M'})(c^f_{M,M''}\otimes\id_{M'})(\id_M\otimes c^f_{M',M''}).
\end{eqnarray*}
\begin{rem}
\label{rem_hexagoncondition}
If the map $f$ satisfies additionally for all $\lambda,\mu,\nu\in\Lambda$
\[ f(\lambda+\nu,\mu)=f(\lambda,\mu)f(\nu,\mu)\:\mbox{and}\:f(\lambda,\mu+\nu)=f(\lambda,\mu)f(\lambda,\nu),\]
then the maps $c^f_{M,M'}$ satisfy the hexagon identities. For more information on the hexagon identities see \cite[XIII 1.1.]{Kassel}.
\end{rem}
\begin{exmp}
There exists a $d\in\N$ such that 
\[ (\lambda,\mu)\in \frac{1}{d}\Z\:\mbox{for all}\: \lambda,\mu\in\Lambda\]
(choose $d$ as the determinant of the Cartan matrix). Assume there is a $v\in k$ such that $v^d=q$. Then for all $w\in k^\times$ the map given by
\[f(\lambda,\mu) := w v^{-d(\lambda,\mu)}\]
satisfies the condition (\ref{equ_f}). If moreover $w=1$, then also the condition from remark \ref{rem_hexagoncondition} is satisfied.
\end{exmp}

\section{Turning $U$-modules into ${}^{\hat{U}}_{\hat{U}}\mYD$-modules (with the same braiding)}
\label{sect_ydstruct}
In this section a Hopf algebra extension $\hat{U}$ of $U_q^{\leq 0}(\mathfrak{g})$ is defined such that every integrable $U$-module is a Yetter-Drinfel'd module over $\hat{U}$ with the property that the induced braiding is the map $c^f$ defined in the section above. A similar construction is mentioned in \cite{Rosso_invent}. Keep the notation from section \ref{sect_uqgmod}, not necessarily assuming that $f$ satisfies the condition from remark \ref{rem_hexagoncondition}. Again let $\operatorname{char} k=0$.\par
It is well known that $U_q^{\leq 0}(\mathfrak{g})$ decomposes as a Radford biproduct
\[ U_q^{\leq 0}(\mathfrak{g}) = \Nichols(V) \# k\Gamma. \]
Here $\Gamma\cong \Z\Phi$ is written multiplicatively identifying $\mu\in\Z\Phi$ with $K_\mu\in\Gamma$ as usual.
$\Nichols(V)$ is the Nichols algebra of the vector space $V := \oplus_{\alpha\in\Pi}k\hat{F}_\alpha$ with braiding
\[ c(\hat{F}_\alpha\otimes \hat{F}_\beta) = q^{-(\alpha,\beta)}\hat{F}_\beta\otimes \hat{F}_\alpha. \]
The usual generators $F_\alpha$ as in Jantzens book are given by $F_\alpha = K_\alpha^{-1} \hat{F}_\alpha$.
The following easy lemma allows to define representations of the biproduct algebra.
\begin{lem}
\label{lem_repbiprod}
Let $H$ be a Hopf algebra with bijective antipode and $R$ a Hopf algebra in ${}_H^H\mYD$. Let $A$ be any algebra. The following data are equivalent:
\begin{itemize}
\item an algebra morphism $\psi:R\# H\rightarrow A$
\item algebra morphisms $\rho:H\rightarrow A$ and $\varphi: R\rightarrow A$ such that:
\[ \forall h\in H,r\in R:\rho(h)\varphi(r) = \varphi(h\s1 \cdot r)\rho(h\s2), \]
where $r$ resp. $h$ run through a set of algebra generators of $R$ resp. $H$.
\end{itemize}
In this case $\psi=\varphi\# \rho$ and $\varphi = \psi|R\# 1,\rho = \psi|1\# H$. 
\end{lem}
\begin{pf}
The proof is straightforward and will be omitted.\qed
\end{pf}
Now the Hopf algebra $\hat{U}$ and the Yetter-Drinfel'd module structure will be constructed in 6 steps.\par
\textbf{Step 1: Enlarge the Group.}\newline
As $\Z\Phi\subset\Lambda$ are free abelian groups of the same rank $|\Pi|$ the quotient $\Lambda/\Z\Phi$ is a finite set. Choose a set $X\subset\Lambda$ of representatives of the cosets of $\Z\Phi$. Define
\[ G := \Gamma \times H, \]
where $H$ denotes the free abelian group generated by the set $X$ (written multiplicatively).
For every $\lambda\in\Lambda$ there are unique elements $\alpha_\lambda\in\Z\Phi$ and $x_\lambda\in X$ such that
\[ \lambda = \alpha_\lambda + x_\lambda.\]
Define for any $\lambda\in\Lambda$
\[ L_\lambda := (K_{\alpha_\lambda}^{-1},x_\lambda) \in G.\]
Note that for $\mu\in\Z\Phi,\lambda\in\Lambda$
\[ L_{\lambda-\mu} = L_\lambda K_\mu.\]

\textbf{Step 2: Define $\hat{U}$}\newline
Now define a $kG$ coaction on $V$ by setting
\[\delta_V(\hat{F}_\alpha) := K_{\alpha}\otimes \hat{F}_\alpha \:\:\mbox{for all}\:\: \alpha\in\Pi.\]
Consider the $kG$-action defined by
\[ K_\alpha \hat{F}_\beta := q^{-(\beta,\alpha)} \hat{F}_\beta \:\:\mbox{and}\:\:L_x \hat{F}_\beta := q^{(\beta,x)} \hat{F}_\beta \]
for all  $\alpha,\beta\in\Pi,x\in X$.
Of course this defines a ${}^G_G\mYD$ structure on $V$ inducing the original braiding. The desired Hopf algebra is
\[ \hat{U} := \Nichols(V)\# kG. \]

\textbf{Step 3: The action of $\hat{U}$ on $U$-modules}\newline
Let $M$ be an integrable $U$-module. Define the action of $G$ on $m\in M_\lambda$ by
\[ K_\alpha m := q^{(\lambda,\alpha)} m \:\:\mbox{and}\:\: L_x m := f(\lambda,x) m \]
for $\alpha\in\Pi,x\in X$. Furthermore consider the action of $\Nichols(V)\subset U$ given by the restriction of the action of $U$ on $M$. Using $\hat{F}_\alpha M_\lambda \subset M_{\lambda - \alpha}$ and the properties of the map $f$ it is easy to check that these two representations satisfy the compatibility conditions from lemma \ref{lem_repbiprod} and induce a representation of $\hat{U}$ on $M$.

\textbf{Step 4: The $\hat{U}$-coaction on $U$-modules}\newline
Let $M$ be an integrable $U$-module. The map
\[ \delta:M\rightarrow \hat{U}\otimes M,\:\: \delta(m) = \sum\limits_{\mu\geq 0} \Theta_\mu^- L_\lambda\otimes \Theta_\mu^+ m\:\:\mbox{for}\:\:m\in M_\lambda\]
defines a coaction on $M$. Of course this map is counital. For $m\in M_\lambda$ calculate
\begin{eqnarray*}
(\id\otimes\delta)\delta(m) &=& \sum\limits_{\nu\geq0} \Theta_\nu^-L_\lambda\otimes\delta(\Theta_\nu^+ m)\\
&=& \sum\limits_{\mu,\nu\geq0} \Theta_\nu^-L_\lambda\otimes\Theta_\mu^-L_\lambda K_\nu^{-1}\otimes\Theta_\mu^+\Theta_\nu^+ m\\
&=& \sum\limits_{\eta\geq 0}\Delta(\Theta_\eta^-L_\lambda)\otimes\Theta_\eta^+ m.
\end{eqnarray*}
In the last step use the equality 
\[ \Delta(\Theta_\eta^-) \otimes \Theta_\eta^+ = \sum\limits_{{\mu,\nu\geq 0}\atop{\mu+\nu=\eta}} \Theta_\nu^-\otimes\Theta_\mu^-K_\nu^{-1}\otimes\Theta_\mu^+\Theta_\nu^+\]
taken from \cite[7.4]{Jantzen}, which holds in $U^{\leq 0}\otimes U^{\leq 0}\otimes U^{>0}\subset \hat{U}\otimes \hat{U}\otimes U^{> 0}$.

\textbf{Step 5: This defines a ${}^{\hat{U}}_{\hat{U}}\mYD$-structure on $M$}\newline
Let $m\in M_\lambda$. It suffices to check the compatibility condition for algebra generators of $\hat{U}$. Start with the $K_\alpha$:
\begin{eqnarray*}
\delta(K_\alpha m) &=& q^{(\lambda,\alpha)}\delta(m)\\
&=& \sum\limits_{\mu\geq0}q^{-(\mu,\alpha)} \Theta_\mu^-L_\lambda\otimes q^{(\lambda+\mu,\alpha)} \Theta_\mu^+ m\\
&=& \sum\limits_{\mu\geq0}K_\alpha \Theta_\mu^-L_\lambda K_\alpha^{-1}\otimes K_\alpha \Theta_\mu^+ m.
\end{eqnarray*}
Then the $L_x$ for $x\in X$:
\begin{eqnarray*}
\delta(L_x m) &=& f(\lambda,x) \delta(m)\\
&=& \sum\limits_{\mu\geq0} q^{(\mu,x)}\Theta_\mu^-L_\lambda\otimes f(\lambda+\mu,x)\Theta_\mu^+m\\
&=& \sum\limits_{\mu\geq0} L_x\Theta_\mu^-L_\lambda L_x^{-1}\otimes L_x\Theta_\mu^+m.\\
\end{eqnarray*}
Finally consider the $F_\alpha,\alpha\in\Pi$.
\[ \delta(F_\alpha m) = \sum\limits_{\mu\geq 0} \Theta_\mu^- L_{\lambda-\alpha}\otimes \Theta_\mu^+F_\alpha m.\]
On the other hand  (setting $\Theta_\mu := 0$ for $\mu\not\geq 0$)
\begin{eqnarray*}
F_\alpha\s1 &m\sm1&S(F_\alpha\s3)\otimes F_\alpha\s2m\sm0 =\\
&=& \phantom{-} \sum\limits_{\mu\geq0}F_\alpha\Theta_\mu^-L_\lambda K_\alpha\otimes K_\alpha^{-1}\Theta_\mu^+m
+\sum\limits_{\mu\geq0}\Theta_\mu^-L_\lambda K_\alpha\otimes F_\alpha\Theta_\mu^+m\\
&&-\sum\limits_{\mu\geq0}\Theta_\mu^-L_\lambda F_\alpha K_\alpha\otimes \Theta_\mu^+ m\\
&=& \phantom{-} \sum\limits_{\mu\geq0}F_\alpha\Theta_{\mu-\alpha}^-L_\lambda K_\alpha\otimes K_\alpha^{-1}\Theta_{\mu-\alpha}^+m
+\sum\limits_{\mu\geq0}\Theta_\mu^-L_\lambda K_\alpha\otimes F_\alpha\Theta_\mu^+m\\
&&-\sum\limits_{\mu\geq0}\Theta_{\mu-\alpha}^-F_\alpha L_\lambda K_\alpha \otimes \Theta_{\mu-\alpha}^+K_\alpha m,\\
\end{eqnarray*}
using $\Delta(F_\alpha) = F_\alpha\otimes K_\alpha^{-1} + 1\otimes F_\alpha$, $S(F_\alpha)= -F_\alpha K_\alpha$ and the commutation relations for the $K_\alpha$'s and $F_\alpha$'s. Now use
\[ \Theta_\mu^-\otimes F_\alpha\Theta_\mu^+ + F_\alpha\Theta_{\mu-\alpha}^-\otimes K_\alpha^{-1}\Theta_{\mu-\alpha}^+ - \Theta_{\mu-\alpha}^-F_\alpha\otimes\Theta_{\mu-\alpha}^+K_\alpha = \Theta_\mu^-\otimes\Theta_\mu^+F_\alpha\]
for all $\mu\geq 0$ from \cite[7.1]{Jantzen}. This yields
\begin{eqnarray*}
F_\alpha\s1 &m\sm1&S(F_\alpha\s3)\otimes F_\alpha\s2m\sm0 =\\
&=& \sum\limits_{\mu\geq0}\Theta_\mu^-L_\lambda K_\alpha\otimes\Theta_\mu^+F_\alpha m \\
&=& \sum\limits_{\mu\geq0}\Theta_\mu^-L_{\lambda-\alpha}\otimes\Theta_\mu^+F_\alpha m = \delta(F_\alpha m).
\end{eqnarray*}

\textbf{Step 6: The induced braiding is $c^f$}\newline
Let $m\in M_\lambda,m'\in M_{\lambda'}$. The braiding induced by the Yetter-Drinfel'd structure defined above is
\begin{eqnarray*}
c_\mYD(m\otimes m') &=& \sum\limits_{\mu\geq0}\Theta_\mu^-L_\lambda m'\otimes\Theta_\mu^+ m\\
 &=& f(\lambda',\lambda)\sum\limits_{\mu\geq0} \Theta_\mu^- m'\otimes\Theta_\mu^+ m\\
&=& c^f(m\otimes m').
\end{eqnarray*}

\begin{rem}
As every $U$-linear map between integrable $U$-modules is $\hat{U}$-linear and colinear (with respect to the structures defined above), this defines a functor from the category of integrable $U$-modules to the category ${}^{\hat{U}}_{\hat{U}}\mYD$. Note that this functor preserves the braiding but is in general \emph{not} monoidal. This is due to the fact that $L_{\lambda+\lambda'} \neq L_\lambda L_{\lambda'}$. In fact if this functor were monoidal, then $c^f$ would satisfy the hexagon identities on every triple of integrable $U$-modules. This is not true unless the condition in remark \ref{rem_hexagoncondition} holds. However if this condition holds there is an other extension $\hat{U}'$ of $U$ and a monoidal functor from the category of integrable $U$-modules to ${}^{\hat{U}'}_{\hat{U}'}\mYD$ that preserves the braiding. In this case choose $G\cong \Lambda$ identifying $\lambda\in\Lambda$ with $K_\lambda\in G$, use $L_\lambda := K_{-\lambda}$ and redo the proof above.
\end{rem}

\begin{rem}
By using similar methods one can find an extension $\hat{U}''$ of $U^{\geq 0}$ and a functor from the category of integrable $U$-modules to ${}^{\hat{U}''}_{\hat{U}''}\mYD$ such that the induced braiding is $(c^f)^{-1}$. Again this functor can not be chosen monoidal unless $f$ has the property from remark \ref{rem_hexagoncondition}.
\end{rem}

\section{Braided biproducts}
\label{sect_brbiprod}
In this section a braided version of Radfords biproduct is introduced. This is done for arbitrary braided categories in \cite{Besp_crossedmodules}. Here an ad-hoc approach for the category ${}^H_H\mYD$ is presented, that leads very quickly to the necessary results. Let $H$ be a Hopf algebra with bijective antipode and $R$ a Hopf algebra in ${}^H_H\mYD$ such that $R\#H$ has bijective antipode. Moreover let $Q$ be a Hopf algebra in ${}^{R\#H}_{R\#H}\mYD$. Consider the projection of Hopf algebras
\[ \eps\otimes\eps\otimes \id_H: Q\#(R\#H)\rightarrow H.\]

\begin{prop}
The space of (right) coinvariants with respect to $\eps\otimes\eps\otimes \id_H$ is $Q\otimes R\otimes 1$.
\end{prop}
\begin{pf}
One inclusion is trivial. So assume there is a coinvariant 
\[T = \sum\limits_{i=1}^r x_i\# r_i\# h_i \in \left(Q\#(R\#H)\right)^{\co\eps\otimes\eps\otimes H}.\]
The $x_i\otimes r_i$ can be chosen linearly independent. Using the formulas for the comultiplication of the Radford biproduct one obtains
\[ T\otimes 1_H = (\id_Q\otimes \id_R\otimes \id_H\otimes\eps\otimes\eps\otimes \id_H)\Delta(T) = \sum\limits_{i=1}^r x_i\otimes r_i\otimes h_i\s1\otimes h_i\s2. \]
This implies $h_i = \eps(h_i)1$ for all $1\leq i\leq r$ and thus $T\in Q\#R\# 1$.
\qed\end{pf}
\begin{defn}
\label{defn_brbiprod}
Thus $Q\otimes R$ inherits the structure of a Hopf algebra in ${}^H_H\mYD$ from the coinvariants. This object is called the \emph{braided biproduct of $Q$ and $R$} and is denoted by $Q\#R$.
\end{defn}
$Q$ is a subalgebra of $Q\#R$ (via the inclusion $x\mapsto x\#1$) and $R$ is a braided Hopf subalgebra of $Q\#R$. \par

Note that $Q\in {}^H_H\mYD$ via the inclusion $H\rightarrow R\#H$ and the projection
\[ \pi_H:R\#H\rightarrow H,\:r\#h \mapsto \eps(r)h.\]
However $Q$ is in general not a braided Hopf algebra in ${}^H_H\mYD$.\par
By construction of $Q\#R$ it is obvious that
\begin{eqnarray*}
 Q\#(R\#H) &\rightarrow& (Q\#R)\#H\\
x\#(r\#h)&\mapsto& (x\#r)\#h
\end{eqnarray*}
is an isomorphism of Hopf algebras.
\subsection{Structure maps}
The following list contains formulas for the structure maps of $Q\#R$. The proofs are left to the reader. For all $x,x'\in Q,r,r'\in R,h\in H$:
\begin{eqnarray*}
(x\#r)(x'\#r')&=& x\left[(r\su1\#r\su2\sm1)\cdot x'\right]\#r\su2\s0r',\\
\Delta_{Q\#R}(x\#r) &=& x\sxu1\#\theta_R(x\sxu2\smx2)\left[\pi_H(x\sxu2\smx1)\cdot r\su1\right]\otimes x\sxu2\sx0\#r\su2,\\
\delta_H(x\#r) &=& \pi_H(x\smx1)r\sm1\otimes x\sx0\#r\s0,\\
h\cdot(x\#r) &=& \left((1\#h\s1)\cdot x\right)\#h\s2\cdot r.\\
\end{eqnarray*}
Note that the action and coaction correspond to the tensor product of Yetter-Drinfel'd modules over $H$.

\subsection{The braided adjoint action}
For any Hopf algebra $R$ in ${}^H_H\mYD$  the braided adjoint action is defined by
\[ \ad_c:R\otimes R\rightarrow R,\:\ad_c(r)(r') := r\su1\left(r\su2\sm1\cdot r'\right)S_R\left(r\su2\s0\right).\]
In the usual Radford biproduct $R\#H$ the following rules are valid:
\begin{eqnarray*}
&\ad&(1\#h)(1\#h')= 1\#\ad(h)(h'),\\
&\ad&(1\#h)(r\#1)= (h\cdot r)\# 1,\\
&\ad&(r\#1)(r'\#1)= \ad_c(r)(r')\#1
\end{eqnarray*}
for all $r,r'\in R,h,h'\in H$.\par
In the braided biproduct $Q\# R$ the corresponding rules
\begin{eqnarray*}
&\ad_c&(1\#r)(1\#r') = 1\#\ad_c(r)(r'),\\
&\ad_c&(1\#r)(x\#1) = \left((r\#1)\cdot x\right)\# 1,\\
&\ad_c&(x\#1)(x'\#1) = \ad_c(x)(x')\#1
\end{eqnarray*}
hold for all $x,x'\in Q,r,r'\in R$.
In the last equation on the right side the ${}^{R\#H}_{R\#H}\mYD$ structure on $Q$ is used to define $\ad_c$.

\section{Graded Yetter-Drinfel'd modules}
\label{sect_gryd}
For this section assume that $A=\oplus_{n\geq0} A(n)$ is a graded Hopf algebra with bijective antipode. Then $H:=A(0)$ is a Hopf algebra with bijective antipode. In this section the notion of a graded Yetter-Drinfel'd modules over $A$ is defined. This class of Yetter-Drinfel'd modules is the natural context for the extension theorem \ref{thm_nicholsext}.
\begin{defn}
\label{defn_gryd}
$M$ is called a \emph{graded Yetter-Drinfel'd module} (over $A$) if $M\in{}^A_A\mYD$ and it has a grading $M=\oplus_{n\geq1}M(n)$ as a vector space such that the action and the coaction are graded maps with respect to the usual grading on tensor products
\[ (A\otimes M)(n) = \sum\limits_{i+j=n} A(i)\otimes M(j).\]
The subspace $M_H := \{m\in M | \delta(m)\in H\otimes M\}$ is called the \emph{space of highest weight vectors} of $M$.\newline
$M$ is said to be \emph{of highest weight} if it is a graded Yetter-Drinfel'd module, $M_H=M(1)$ and $M$ is generated by $M_H$ as an $A$-module.
\end{defn}

\begin{lem}
Let $M\in {}^A_A\mYD$ be of highest weight. The space of highest weight vectors $M_H$ of $M$ is a Yetter-Drinfel'd module over $H$ with action and coaction given by the restrictions of the structure maps on $M$. 
\end{lem}

\begin{pf}
$M_H$ is an $H$ submodule by the Yetter-Drinfel'd condition. To see that $M_H$ is a $H$-comodule fix a basis $(h_i)_{i\in I}$ of $H$. There are scalars $(\alpha_{jl}^i)_{i,j,l\in I}$ such that for all $i\in I$
\[ \Delta(h_i) = \sum\limits_{j,l\in I} \alpha_{jl}^i h_j\otimes h_l.\]
Furthermore let $m\in M_H$. Now there are elements $(m_i)_{i\in I}$ of $M$ (almost all equal to zero) such that
\[ \delta(m) = \sum\limits_{i\in I} h_i\otimes m_i.\]
It suffices to show that $m_j\in M_H$ for all $j\in I$. It is
\[ \sum\limits_{j\in I} h_j\otimes\delta(m_j) = (H\otimes\delta)\delta(m) = (\Delta\otimes M)\delta(m) = \sum\limits_{i,j,l\in I} \alpha_{jl}^i h_j\otimes h_l\otimes m_i\]
and thus for all $j\in I$
\[ \delta(m_j) = \sum\limits_{i,l\in I} \alpha_{jl}^i  h_l\otimes m_i\in H\otimes M,\]
showing $m_j\in M_H$ for all $j\in I$.
\qed\end{pf}

\begin{exmp}
\label{exmp_uqgmod}
Let $\hat{U}$ be the extension of $U_q^{\leq 0}(\mathfrak{g})$ defined in section \ref{sect_ydstruct} and $M$ a simple integrable $U_q(\mathfrak{g})$-module of highest weight $\lambda$. Define a grading on $M$ by
\[ M(n) := \sum\limits_{{\mu\geq 0}\atop{\operatorname{ht}\mu=n-1}} M_{\lambda - \mu}.\]
The Yetter-Drinfel'd module structure defined in section \ref{sect_ydstruct} makes $M$ a graded Yetter-Drinfel'd module over $\hat{U}$ of highest weight, where $\hat{U}$ inherits its grading from the grading on $\Nichols(V)$. As each finite-dimensional integrable $U_q(\mathfrak{g})$-module is a direct sum of modules of highest weight and because the functor from section \ref{sect_ydstruct} preserves direct sums, one obtains that any finite-dimensional $U_q(\mathfrak{g})$-module is a graded Yetter-Drinfel'd module over $\hat{U}$ of highest weight. The space of highest weight vectors is exactly the space spanned by the vectors that are of highest weight in the usual sense.
\end{exmp}
The next step is to extend the grading from graded Yetter-Drinfel'd modules to their Nichols algebras. In general the coradical grading of the Nichols algebra does not turn the action and coaction into graded maps.
\begin{prop}
\label{prop_grading}
Let $M\in{}^A_A\mYD$ be a graded Yetter-Drinfel'd module. Then there is a grading
\[ \Nichols(M) = \bigoplus\limits_{n\geq 0} \Nichols(M)[n], \]
turning $\Nichols(M)$ into a graded Yetter-Drinfel'd module over $A$ and into a graded braided Hopf algebra such that
\[ \Nichols(M)[0]=k1 \:\:\mbox{and}\:\: \Nichols(M)[1] = M(1).\]
\end{prop}
\begin{pf}
Grade the tensor algebra $T(M)$ by giving $M(n)$ the degree $n$. Then the action and the coaction are graded and so is the braiding. Thus the quantum antisymmetrizer maps are graded maps. As the kernel of the projection $T(M)\rightarrow \Nichols(M)$ is just the direct sum of the kernels of the quantum antisymmetrizers \cite{Sbg_borel}, it is a graded Hopf ideal. Thus the quotient $\Nichols(M)$ admits the desired (induced) grading.
\qed\end{pf}

\section{Braided biproducts of Nichols algebras}
\label{sect_extthm}
In this section the results of the preceeding sections are specialized to a braided biproduct of two Nichols algebras. The next theorem is a generalization of \cite[Proposition 2.2]{Rosso_invent} from abelian group algebras to arbitrary Hopf algebras $H$ with bijective antipode. We give a different proof, using the grading instead of the bilinear form on the Nichols algebra. This result allows to reduce the study of Nichols algebras of graded Yetter-Drinfel'd modules over $\Nichols(V)\#H$ to the study of Yetter-Drinfel'd modules over $H$.
\begin{thm}
\label{thm_nicholsext}
Assume that $H$ is a Hopf algebra with bijective antipode, $V\in {}^H_H\mYD$ and set $A:=\Nichols(V)\# H$ as a graded Hopf algebra with grading $A(n) := \Nichols(V)(n)\# H$. Furthermore let $M\in {}^A_A\mYD$ be a graded Yetter-Drinfel'd module.
If $M$ is of highest weight then
\[ \Nichols(M)\#\Nichols(V) \cong \Nichols(M_H\oplus V)\]
as graded braided Hopf algebras in ${}^H_H\mYD$. Here the left side is graded by the tensor product grading and the grading for $\Nichols(M)$ is taken from proposition \ref{prop_grading}.
\end{thm}
\begin{pf}
$B := \Nichols(M)\#\Nichols(V)$ is graded as a braided Hopf algebra: All the structure maps of $B$ are obtained from the structure maps of $H, V, \Nichols(V), M$ and $\Nichols(M)$. As all these maps are graded (giving $V$ the degree 1 and $H$ the degree 0) this part is done.\par
Next check that $B[1] = P(B)$. As $B$ is graded as a coalgebra and $B[0]=k1$ it is clear that $B[1]\subset P(B)$. To show the other inclusion identify $B$ with the coinvariant subalgebra in $\Nichols(M)\#(\Nichols(V)\# H)$. By construction of the grading 
\[ B[1] = M_H\#1\#1\oplus 1\#V\#1.\]
For a primitive element $t = \sum_{i=0}^r x_i\#u_i\#1 \in P(B)$, the coproduct of $t$ is given by
\[ \Delta(t) = \sum_{i=0}^r x_i\sxu1\#x_i\sxu2\smx2\left(u_i\su1\# 1\right)\iota S\pi\left(x_i\sxu2\smx1\right)\otimes x_i\sxu2\sx0\#u_i\su2\#1,\]
where $\pi$ denotes the Hopf algebra projection from $A=\Nichols(V)\#H$ onto $H$ and $\iota$ is the inclusion of $H$ into $\Nichols(V)\# H$.\newline
The $x_i$ can be chosen linearly independent and such that $\eps(x_i)\neq 0$ if and only if $i=0$. Applying the map $\eps_{\Nichols(M)}$ to the fourth tensorand of the equality $1\otimes t + t\otimes 1 = \Delta(t)$ yields
\[ u_0\in P(\Nichols(V))=V \:\:\mbox{and}\:\forall 1\leq i\leq r: u_i\in k1. \]
This means $t = x\#1\#1 + 1\#u\#1$ for $x\in \Nichols(M),u\in V$. In particular $x\#1\#1\in P(B)$. Now calculate
\[\Delta(x\#1\#1) = x\sxu1\#x\sxu2\smx2\iota S\pi\left(x\sxu2\smx1\right)\otimes x\sxu2\sx0\# 1\#1,\]
which holds in $\Nichols(M)\otimes\Nichols(V)\otimes H\otimes\Nichols(M)\otimes\Nichols(V)\otimes H$. Applying to counits to the second and third tensorand yields $x\in P(\Nichols(M))$. Then apply $\eps_{\Nichols(M)}$ to the first tensorand of the equation and observe that
\[ x\smx2\iota S\pi\left(x\smx1\right)\otimes x\sx0 = 1\otimes x,\]
implying that $\delta(x) = \iota\pi(x\smx1)\otimes x\sx0\in H\otimes M$ and thus $x\in M_H$.\par
It remains to show that $B$ is actually generated by $B[1]$. Of course $B$ is generated by $\Nichols(M)\#1\#1$ and $1\#\Nichols(V)\#1$. So it suffices to show that $M\#1\#1$ is contained in the subalgebra generated by $M_H\#1\#1$ and $1\#V\#1$. As $M$ is of highest weight it is generated as a $\Nichols(V)\# H$-module by $M_H$. Using that $M_H$ is an $H$-module this means
\begin{eqnarray*}
M &=& (\Nichols(V)\#H)\cdot M_H = ((\Nichols(V)\#1)(1\# H))\cdot M_H  \\&=& (\Nichols(V)\#1)\cdot M_H = \ad_c(\Nichols(V))(M_H)\#1.
\end{eqnarray*}
Thus within $B$, $M$ is generated by $M_H$ under the braided adjoint action of $\Nichols(V)$. Alltogether $M_H$ and $V$ generate $B$.
\qed\end{pf}

\section{The Gelfand-Kirillov dimension of Nichols algebras of integrable $U_q(\mathfrak{g})$-modules}
\label{sect_uqgapp}
Assume that $\operatorname{char} k=0$. For this section let $M$ be a finite-dimensional integrable $U_q(\mathfrak{g})$-module with braiding $c^f$ as in section \ref{sect_uqgmod}. The first result will be a criterion to decide whether $\Nichols(M)$ has finite Gelfand-Kirillov dimension (recall that $q$ is not a root of unity). From now on restrict to braidings of the following special form. This restriction is necessary due to missing information on Nichols algebras of diagonal type.

\begin{defn}
The braiding $c^f$ is \emph{of exponential type with function $\varphi$} if the map $f:\Lambda\times\Lambda \rightarrow k^\times$ is of the form
\[ f(\lambda,\mu) = v^{-d\varphi(\lambda,\mu)} \]
for some $v\in k$, $d\in\Z$ such that $v^d=q$ and for a map $\varphi:\Lambda\otimes\Lambda\rightarrow \frac{2}{d}\Z$ such that for $\lambda,\lambda'\in\Lambda,\nu\in\Z\Phi$
\[ \varphi(\lambda+\nu,\lambda') = \varphi(\lambda,\lambda')+(\nu,\lambda')\:\:\mbox{and}\:\:\varphi(\lambda,\lambda'+\nu) = \varphi(\lambda,\lambda')+(\lambda,\nu).\]
\\
$c^f$ is \emph{of strong exponential type} if it is of exponential type with a function $\varphi$ and for all $\lambda,\lambda'\in\Lambda$ that are highest weights of $M$
\[ \varphi(\lambda,\lambda) \leq 0 \: \Rightarrow \: \left(\lambda = 0 \:\:\mbox{and}\:\:\varphi(\lambda,\lambda')+\varphi(\lambda',\lambda)=0\right).\]
\end{defn}

As shown in example \ref{exmp_uqgmod} the module $M$ is a graded Yetter-Drinfel'd module of highest weight over $\hat{U}$ (the grading on $\hat{U} = \Nichols(V)\#kG$ is the one induced by the Nichols algebra). Assume that $M=\oplus_{1\leq i\leq r} M_i$ is the decomposition of the $U_q(\mathfrak{g})$-module $M$ into irreducible submodules. For all $1\leq i\leq r$ choose a highest weight vector $m_i\in M_i$ and denote by $\lambda_i\in\Lambda$ the weight of $m_i$. Then the space $M_{kG}\oplus V$ has basis $\{\hat{F}_\alpha|\alpha\in\Pi\}\cup\{m_1,\ldots,m_r\}$. If $c^f$ is of exponential type, the braiding on $M_{kG}\oplus V$ is given by
\begin{eqnarray*}
c(\hat{F}_\alpha\otimes \hat{F}_\beta) &=& v^{-d(\beta,\alpha)} \hat{F}_\beta\otimes \hat{F}_\alpha,\\
c(\hat{F}_\alpha\otimes m_j) &=& v^{d(\lambda_j,\alpha)} m_j\otimes \hat{F}_\alpha,\\
c(m_i\otimes \hat{F}_\beta) &=& v^{d(\beta,\lambda_i)} \hat{F}_\beta\otimes m_i\:\:\mbox{and}\\
c(m_i\otimes m_j) &=& f(\lambda_j,\lambda_i) m_j\otimes m_i = v^{-d\varphi(\lambda_i,\lambda_j)}m_j\otimes m_i.
\end{eqnarray*}

Let $P := \Pi\dot{\cup}\{1,\ldots,r\}$. If $c^f$ is of strong exponential type there is always a matrix $(b_{ij})_{i,j\in P}\in \Q^{P\times P}$ such that the following conditions are satisfied:
\begin{eqnarray*}
\forall\alpha,\beta\in\Pi:& 2(\alpha,\beta) = (\alpha,\alpha)b_{\alpha\beta}&\:\:(1)\\
\forall\alpha\in\Pi,1\leq i\leq r:& 2(\alpha,\lambda_i) = -(\alpha,\alpha)b_{\alpha i}&\:\:(2)\\
\forall\alpha\in\Pi,1\leq i\leq r:& 2(\alpha,\lambda_i) = -\varphi(\lambda_i,\lambda_i)b_{i\alpha}&\:\:(3)\\
\forall 1\leq i,j\leq r:& \varphi(\lambda_i,\lambda_j)+\varphi(\lambda_j,\lambda_i) = \varphi(\lambda_i,\lambda_i)b_{ij}&\:\:(4)\\
\forall 1\leq i\leq r,i\neq j\in P,\alpha\in\Pi:& \varphi(\lambda_i,\lambda_i)=0\Rightarrow b_{ii}=2, b_{ij}=0, b_{i\alpha}=0&\:\:(5)
\end{eqnarray*}
The matrix $(b_{ij})_{i,j\in P}$ will be called the \emph{extended Cartan matrix} of $M$.

\begin{thm}
\label{thm_finkrit}
Assume that the braiding $c^f$ on the finite-dimensional integrable $U_q(\mathfrak{g})$-module $M$ is of exponential type with a symmetric function $\varphi$ (i.e. $\varphi(\lambda,\lambda') = \varphi(\lambda',\lambda)$ for all $\lambda,\lambda'\in\Lambda$). \newline
 Then the Nichols algebra $\Nichols(M,c^f)$ has finite Gelfand-Kirillov dimension if and only if $c^f$ is of strong exponential type (with function $\varphi$) and the extended Cartan matrix $(b_{ij})$ is a Cartan matrix of finite type.
\end{thm}
\begin{pf}
Denote the basis of $M_{kG}\oplus V$ by $x_i,i\in P$, where $x_\alpha := \hat{F}_\alpha$ and $x_i := m_i$ for $\alpha\in\Pi$ and $1\leq i\leq r$. The braiding of $M_{kG}\oplus V$ is of the form
\[c(x_i\otimes x_j) = q_{ij} x_j\otimes x_i \:\:\forall i,j\in P,\]
where the $q_{ij}$ can be read off the formulas given above:
\[ q_{\alpha\beta} = v^{-d(\beta,\alpha)}, q_{\alpha i} = v^{d(\alpha,\lambda_i)}, q_{i\alpha} = v^{d(\alpha,i)}\:\mbox{and}\: q_{ij} = v^{-d\varphi(\lambda_i,\lambda_j)}\]
for all $\alpha,\beta\in\Pi$ and $1\leq i,j\leq r$.\par
\textbf{The if-part:}  By the definition of $(b_{ij})$ for all $i,j\in P$
\[ q_{ij}q_{ji}= {q_{ii}}^{b_{ij}}.\]
For all $\alpha\in \Pi$ and for all $1\leq i\leq r$ with $\varphi(\lambda_i,\lambda_i)\neq 0$ define 
\[ d_\alpha := \frac{d(\alpha,\alpha)}{2}\:\:\mbox{and}\:\: d_i := \frac{d\varphi(\lambda_i,\lambda_i)}{2}\]
and for $1\leq i\leq r$ with $\varphi(\lambda_i,\lambda_i)=0$ define $d_i := 1$.
These $(d_i)$ are positive integers satisfying
\[ d_ib_{ij} = d_jb_{ji}\:\:\forall i,j\in P.\]
Because $\varphi$ is symmetric one gets for all $i,j\in P$
\[ q_{ij} = v^{-d_ib_{ij}}.\]
This means that the braiding on $M_{kG}\oplus V$ is of Drinfeld-Jimbo type as defined in \cite{AS4} with generalized Cartan matrix $(b_{ij})$. As $(b_{ij})$ is a finite Cartan matrix, $\Nichols(M_{kG}\oplus V)$ has finite Gelfand-Kirillov dimension by \cite[Theorem 2.10.]{AS4}.
Thus the subalgebra $\Nichols(M)$ has finite Gelfand-Kirillov dimension \cite[Lemma 3.1.]{KL}.\par
\textbf{The only-if-part:} By \cite[Theorem 36]{Ich_PBW} $\Nichols(M)$ has a PBW basis and because the Gelfand-Kirillov dimension is finite the set of PBW generators $P_M$ must be finite. Similarily $\Nichols(V)$ has a PBW basis and because it has finite Gelfand-Kirillov dimension (see \cite[Theorem 2.10.]{AS4}) its set of PBW generators $P_V$ must be finite. Thus the set 
\[ P := \{p\#1|p\in P_M\}\dot{\cup}\{1\# p'|p'\in P_V\}\]
forming a set of PBW generators for $\Nichols(M)\#\Nichols(V) \cong \Nichols(M_{kG}\oplus V)$ is finite. Hence $\Nichols(M_{kG}\oplus V)$ has finite Gelfand-Kirillov dimension. Now \cite[Lemma 14 and 20]{Rosso_invent} allow to find integers $c_{ij}\leq 0, {i,j\in P}$ such that
\[ q_{ij}q_{ji} = q_{ii}^{c_{ij}}\:\:\forall i,j\in P.\]
Using the definition of the $q_{ij}$ one obtains that the $c_{ij}$ must satisfy the equations $(1)-(4)$ from the definition of $(b_{ij})$ with $b_{ij},i,j\in P$ replaced by $c_{ij},i,j\in P$ ($v$ is not a root of unity). Because of relations $(3)$ and $(4)$, $\varphi$ must satisfy the condition from the definition of strong exponential braidings. Furthermore one may assume $c_{ii}=2$ for all $i\in P$ and $c_{ij}=0$ for all $1\leq i\leq r$ with $\varphi(\lambda_i,\lambda_i)=0, i\neq j\in P$. This means that $b_{ij}=c_{ij}$ for all $i,j\in P$. Now observe that $b_{ij}$ is a generalized Cartan matrix. Exactly as in the ``only-if'' part of the proof the braiding in $M_{kG}\oplus V$ is of Drinfel'd-Jimbo type with generalized Cartan matrix $(b_{ij})$. By \cite[Theorem 2.10.]{AS4} $(b_{ij})$ is a finite Cartan matrix because $\Nichols(M_{kG}\oplus V)$ has finite Gelfand-Kirillov dimension.
\qed\end{pf}

\subsection{Explicit calculations for finite-dimensional simple $U_q(\mathfrak{g})$-modules}
In this section the results above are used to determine all pairs $(\lambda,\varphi)$ for each finite-dimensional simple complex Lie algebra $\mathfrak{g}$ such that the Nichols algebra of the $U_q(\mathfrak{g})$-module of highest weight $\lambda$ together with the braiding defined by the function $\varphi$ has finite Gelfand-Kirillov dimension. First observe that (as only modules of highest weight are considered) one may assume that the function $\varphi$ is of the form
\[ \varphi(\mu,\nu) = (\mu,\nu)+x\:\:\mbox{for}\:\:\mu,\nu\in\Lambda\]
for some $x\in\frac{1}{d}\Z$, where $d$ is the determinant of the Cartan matrix. This is true because the braiding $c^f$ depends only on the values $\varphi(\lambda',\lambda'')$ for those weights $\lambda',\lambda''\in\Lambda$ such that $M_{\lambda'}\neq 0,M_{\lambda''}\neq 0$. Thus one can choose
\[ x:= \varphi(\lambda,\lambda)-(\lambda,\lambda)\]
for any weight $\lambda\in\Lambda$ with $M_\lambda\neq 0$.\par

\begin{thm}
Let $\mathfrak{g}$ be a finite-dimensional simple complex Lie algebra with weight lattice $\Lambda$. Fix a $U_q(\mathfrak{g})$-module $M$ of highest weight $\lambda\in\Lambda$ and a value $x\in \Q$. Let $d'$ be the least common multiple the denominator of $x$ and the determinant of the Cartan matrix of $\mathfrak{g}$. Let $d:=2d'$ and fix $v\in k$ with $v^d = q$. Define a function
\[ f:\Lambda\times\Lambda\rightarrow k^\times,\:\: (\lambda,\lambda') \mapsto v^{d\left((\lambda,\lambda') + x\right)}.\]
The Nichols algebra $\Nichols(M,c^f_{M,M})$ has finite Gelfand-Kirillov dimension if and only if the tuple $\mathfrak{g},\lambda,x$ occurs in table \ref{table_extensions}.
\end{thm}

\begin{table}
\begin{center}
\begin{tabular}{c|c|c|c|c|c}
Type of $\mathfrak{g}$ & $\lambda$ & $x$ & Type of $(b_{ij})$ & $\varphi(\lambda,\lambda)$ & relations in degree\\
\hline
\hline
any $D$ & 0 & any & $D\cup A_0$ & any & no relations\\ \hline
$A_n, n\geq 1$ & $\lambda_1, \lambda_n$ & $\frac{n+2}{n+1}$ & $A_{n+1}$ & 2 & 2\\\hline
$A_n, n\geq 1$ & $\lambda_1, \lambda_n$ & $\frac{1}{n+1}$ & $B_{n+1}$ & 1 & 2\\\hline
$A_n, n\geq 1$ & $2\lambda_1, 2\lambda_n$ & $\frac{4}{n+1}$ & $C_{n+1}$ (resp. $B_2$) & 4 & 3\\\hline
$A_n, n\geq 3$ & $\lambda_{n-1},\lambda_2$ & $\frac{4}{n+1}$ & $D_{n+1}$ & 2 & 2\\\hline
$A_1$ & $\lambda_1$ & $\frac{1}{6}$ & $G_2$ & $\frac{3}{2}$ & 2\\\hline
$A_1$ & $3\lambda_1$ & $\frac{3}{2}$ & $G_2$ & 6 & 4\\\hline
$A_5$ & $\lambda_3$ & $\frac{1}{2}$ & $E_6$ & 2 & 2\\\hline
$A_6$ & $\lambda_3, \lambda_4$ & $\frac{5}{7}$ & $E_7$&2 & 2\\\hline
$A_7$ & $\lambda_3, \lambda_5$ & $\frac{1}{8}$ & $E_8$&2 & 2\\\hline
$B_n, n\geq 2$ & $\lambda_1$ & $1$ & $B_{n+1}$ & 2 & 2\\\hline
$C_n, n\geq 3$ & $\lambda_1$ & $1$ & $C_{n+1}$ & 2 & 2\\\hline
$C_3$ & $2\lambda_3$ & $-2$ & $F_4$ & 4 & 2\\\hline
$D_4$ & $\lambda_1, \lambda_3, \lambda_4$ & $1$ & $D_5$ & 2 & 2\\\hline
$D_n ,n\geq 5$ & $\lambda_1$ & $1$ & $D_{n+1}$ & 2 & 2\\\hline
$D_5$ & $\lambda_4, \lambda_5$ & $\frac{3}{4}$ & $E_6$ & 2 & 2\\\hline
$D_6$ & $\lambda_5, \lambda_6$ & $\frac{1}{2}$ & $E_7$ & 2 & 2\\\hline
$D_7$ & $\lambda_6, \lambda_7$ & $\frac{1}{4}$ & $E_8$ & 2 & 2\\\hline
$E_6$ & $\lambda_1, \lambda_6$ & $\frac{2}{3}$ & $E_7$ & 2 & 2\\\hline
$E_7$ & $\lambda_7 $ & $1$ & $E_8$ & 2 & 2\\
\end{tabular}
\caption{\label{table_extensions}
Highest weights with Nichols algebras of finite Gelfand-Kirillov dimension}\end{center}
\end{table}
Note that the braiding $c^f_{M,M}$ does depend on the choice of $v$, but the Gelfand-Kirillov dimension of $\Nichols(M,c^f_{M,M})$ does not.\newline
In table \ref{table_extensions} also the type of the extended Cartan matrix $(b_{ij})$ and the value $\varphi(\lambda,\lambda)$ are given. The weight $\lambda_i$ always denotes the fundamental weight dual to the root $\alpha_i$. The numbering of the roots is as in \cite{Humphreys}.
\begin{pf}
Assume that the tuple $\mathfrak{g},\lambda,x$ leads to a Nichols algebra of finite Gelfand-Kirillov dimension and let $(a_{\alpha\beta})_{\alpha,\beta\in\Pi}$ be the Cartan matrix for $\mathfrak{g}$. By theorem \ref{thm_finkrit} the extended Cartan matrix $(b_{ij})_{i,j\in P}$ is a finite Cartan matrix. Furthermore $P=\Pi\dot{\cup}\{1\}$ and $b_{\alpha\beta} = a_{\alpha\beta}$ for $\alpha,\beta\in\Pi$. First assume that $(b_{ij})$ is not a connected Cartan matrix. As $(a_{ij})$ is a connected Cartan matrix observe
\[ b_{1\alpha} = 0 = b_{\alpha 1}\:\:\forall\alpha\in\Pi.\]
By the definition of $(b_{ij})$ this imples $\lambda=0$. This is the first line in the table.\par
Now assume that $(b_{ij})$ is a connected finite Cartan matrix and thus its Coxeter graph contains no cycles. As $(a_{\alpha\beta})$ is also a connected finite Cartan matrix there is a unique root $\alpha\in\Pi$ such that
\[ b_{\alpha 1},b_{1\alpha}< 0 \:\:\mbox{and for all}\:\beta\in\Pi\setminus\{\alpha\}: b_{\beta 1}=0=b_{1 \beta}.\]
This implies that
\[ l:=(\alpha,\lambda)> 0\:\:\mbox{and for all}\:\beta\in\Pi\setminus\{\alpha\}: (\beta,\lambda)=0.\]
Conclude that $\lambda = l\lambda_\alpha$ and $l\in\N$, where $\lambda_\alpha$ is the weight dual to the root $\alpha$, i.e.
\[ (\lambda_\alpha,\beta) = \delta_{\beta,\alpha}.\]
Now observe, using the definition of $(b_{ij})$ and $\varphi$, that
\[ l= -\frac{b_{\alpha 1}(\alpha,\alpha)}{2},\]
\[ \varphi(\lambda,\lambda) = \frac{b_{\alpha 1}}{b_{1 \alpha}}(\alpha,\alpha),\:\:\mbox{and}\]
\[ x = \varphi(\lambda,\lambda) - (\lambda,\lambda).\]
Now in a case-by-case analysis consider all finite connected Cartan matrices $(a_{\alpha\beta})_{\alpha,\beta\in\Pi}$ and all possible finite connected Cartan matrices $(b_{ij})_{i,j\in P}$ having $(a_{\alpha\beta})$ as a submatrix. In each case compute the values for $l,\varphi(\lambda,\lambda)$ and $x$ and decide if there is a tuple $\mathfrak{g},\lambda,x$ leading to the matrix $(b_{ij})$. For every case also the Dynkin diagram of $(b_{ij})$ with labeled vertices is given. The vertices $1,\ldots,n$ correspond to the simple roots $\alpha_1,\ldots,\alpha_n\in\Pi$, the vertex $\star$ corresponds to $1\in P$.\par
\underline{$A_n\rightarrow A_{n+1},n\geq 1$}:
\begin{diagram}[size=1em,bottom]
\stackrel{1}{\bullet}&\rLine&\stackrel{2}{\bullet}&\cdots\cdots&\stackrel{n-1}{\bullet}&\rLine&\stackrel{n}{\bullet}&\rLine&\stackrel{\star}{\bullet}&
\:\:\mbox{or}\:\:&
\stackrel{\star}{\bullet}&\rLine&\stackrel{1}{\bullet}&\rLine&\stackrel{2}{\bullet}&\cdots\cdots&\stackrel{n-1}{\bullet}&\rLine&\stackrel{n}{\bullet}\\
\end{diagram}
It is either $\alpha = \alpha_n$ or $\alpha = \alpha_1$. In any case 
\[ l = 1,\varphi(\lambda,\lambda) = 2, x = \frac{n+2}{n+1}.\]
So $\lambda = \lambda_1$ or $\lambda = \lambda_n$ together with $ x = \frac{n+2}{n+1}$ extend the type $A_n$ to $A_{n+1}$.

\underline{$A_n\rightarrow B_{n+1},n\geq 1$}:
\begin{diagram}[size=1em,bottom]
\stackrel{1}{\bullet}&\rLine&\stackrel{2}{\bullet}&\cdots\cdots&\stackrel{n-1}{\bullet}&\rLine&\stackrel{n}{\bullet}&\rImplies&\stackrel{\star}{\bullet}&
\:\:\mbox{or}\:\:&
\stackrel{\star}{\bullet}&\lImplies&\stackrel{1}{\bullet}&\rLine&\stackrel{2}{\bullet}&\cdots\cdots&\stackrel{n-1}{\bullet}&\rLine&\stackrel{n}{\bullet}\\
\end{diagram}
It is either  $\alpha = \alpha_1$ or $\alpha = \alpha_n$ and obtain
\[ l = 1,\varphi(\lambda,\lambda) = 1, x = \frac{1}{n+1}.\]
So $\lambda = \lambda_1$ or $\lambda = \lambda_n$ together with $ x = \frac{1}{n+1}$ extend the type $A_n$ to $B_{n+1}$.

\underline{$A_n \rightarrow C_{n+1}, n\geq 1$}:
\begin{diagram}[size=1em,bottom]
\stackrel{1}{\bullet}&\rLine&\stackrel{2}{\bullet}&\cdots\cdots&\stackrel{n-1}{\bullet}&\rLine&\stackrel{n}{\bullet}&\lImplies&\stackrel{\star}{\bullet}&
\:\:\mbox{or}\:\:&
\stackrel{\star}{\bullet}&\rImplies&\stackrel{1}{\bullet}&\rLine&\stackrel{2}{\bullet}&\cdots\cdots&\stackrel{n-1}{\bullet}&\rLine&\stackrel{n}{\bullet}\\
\end{diagram}
It is either  $\alpha = \alpha_1$ or $\alpha = \alpha_n$ and calculate
\[ l = 2, \varphi(\lambda,\lambda) = 4, x = \frac{4}{n+1}.\]
So $\lambda = 2\lambda_1$ or $\lambda = 2\lambda_n$ together with $ x = \frac{4}{n+1}$ extend the type $A_n$ to $C_{n+1}$ (resp. $B_2$).

\underline{$A_n\rightarrow D_{n+1}, n\geq 3$}:
\begin{diagram}[size=1em,bottom]
\stackrel{1}{\bullet}&\rLine&\stackrel{2}{\bullet}&\cdots\cdots&\stackrel{n-1}{\bullet}&\rLine&\stackrel{n}{\bullet}&
\:\:\mbox{or}\:\:&
\stackrel{1}{\bullet}&\rLine&\stackrel{2}{\bullet}&\cdots\cdots&\stackrel{n-1}{\bullet}&\rLine&\stackrel{n}{\bullet}\\
&&&&\dLine&&& &&&\dLine&&&&\\
&&&&\bullet&&& &&&\bullet&&&&\\
&&&&{}^\star&&& &&&{}^\star&&&&
\end{diagram}
It is either  $\alpha = \alpha_2$ or $\alpha = \alpha_{n-1}$ and calculate
\[ l = 1, \varphi(\lambda,\lambda) = 2, x = \frac{4}{n+1}.\]
So $\lambda = \lambda_2$ or $\lambda = \lambda_{n-1}$ together with $ x = \frac{4}{n+1}$ extend the type $A_n$ to $D_{n+1}$.

\underline{$A_1\rightarrow G_2$}:
\begin{diagram}[size=1em,bottom]
\stackrel{1}{\bullet}&\rThree&\stackrel{\star}{\bullet}&\:\:\mbox{or}\:\:&
\stackrel{1}{\bullet}&\lThree&\stackrel{\star}{\bullet}
\end{diagram}
In any case $\alpha = \alpha_1$. The left diagram leads to
\[ l = 1, \varphi(\lambda,\lambda) = \frac{2}{3}, x = \frac{1}{6}.\]
The right diagram corresponds to
\[ l = 3, \varphi(\lambda,\lambda) = 6, x = \frac{3}{2}.\]
So $\lambda = \lambda_1, x = \frac{1}{6}$ or $\lambda = 3\lambda_{1},x=\frac{3}{2}$ extend the type $A_1$ to $G_2$.

\underline{$A_5\rightarrow E_6$}:
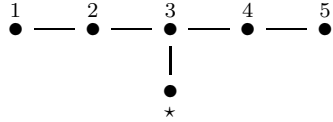
\begin{diagram}[size=1em,bottom]
\stackrel{1}{\bullet}&\rLine&\stackrel{2}{\bullet}&\rLine&\stackrel{3}{\bullet}&\rLine&\stackrel{4}{\bullet}&\rLine&\stackrel{5}{\bullet}\\
&&&&\dLine&&&&\\
&&&&\bullet&&&&\\
&&&&{}^\star&&&&
\end{diagram}
Here $\alpha = \alpha_3$ and thus
\[ l = 1, \varphi(\lambda,\lambda) = 2, x = \frac{1}{2}.\]

\underline{$A_6\rightarrow E_7$}:
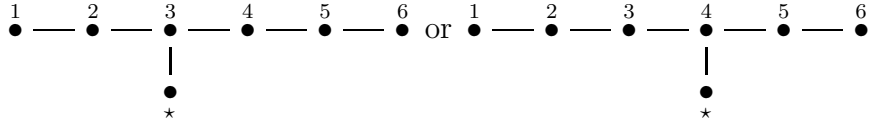
\begin{diagram}[size=1em,bottom]
\stackrel{1}{\bullet}&\rLine&\stackrel{2}{\bullet}&\rLine&\stackrel{3}{\bullet}&\rLine&\stackrel{4}{\bullet}&\rLine&\stackrel{5}{\bullet}&\rLine&\stackrel{6}{\bullet}&
\:\:\mbox{or}\:\:&
\stackrel{1}{\bullet}&\rLine&\stackrel{2}{\bullet}&\rLine&\stackrel{3}{\bullet}&\rLine&\stackrel{4}{\bullet}&\rLine&\stackrel{5}{\bullet}&\rLine&\stackrel{6}{\bullet}\\
&&&&\dLine&&&&&& &&&&&&&&\dLine &&&& \\
&&&&\bullet&&&&&& &&&&&&&&\bullet &&&&\\
&&&&{}^\star&&&&&& &&&&&&&&{}^\star &&&&
\end{diagram}
$\alpha = \alpha_3$ or $\alpha=\alpha_4$ leads to
\[ l = 1, \varphi(\lambda,\lambda) = 2, x = \frac{5}{9}.\]

\underline{$A_7\rightarrow E_8$}:
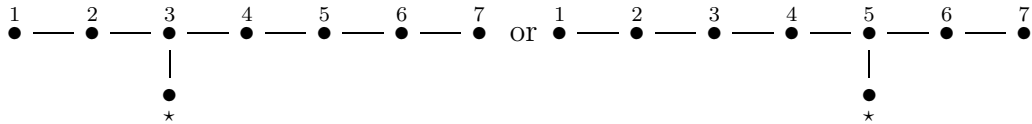
\begin{diagram}[size=1em,bottom]
\stackrel{1}{\bullet}&\rLine&\stackrel{2}{\bullet}&\rLine&\stackrel{3}{\bullet}&\rLine&\stackrel{4}{\bullet}&\rLine&\stackrel{5}{\bullet}&\rLine&\stackrel{6}{\bullet}&\rLine&\stackrel{7}{\bullet}
\:\:\mbox{or}\:\:&
\stackrel{1}{\bullet}&\rLine&\stackrel{2}{\bullet}&\rLine&\stackrel{3}{\bullet}&\rLine&\stackrel{4}{\bullet}&\rLine&\stackrel{5}{\bullet}&\rLine&\stackrel{6}{\bullet}&\rLine&\stackrel{7}{\bullet}\\
&&&&\dLine&&&&&&&& &&&&&&&&&\dLine &&&&&& \\
&&&&\bullet&&&&&&&& &&&&&&&&&\bullet &&&&&&\\
&&&&{}^\star&&&&&&&& &&&&&&&&&{}^\star &&&&&&
\end{diagram}
$\alpha = \alpha_3$ or $\alpha=\alpha_5$ yield
\[ l = 1, \varphi(\lambda,\lambda) = 2, x = \frac{1}{8}.\]

\underline{Almost all of the other cases} follow the same principle and will be omitted. There are two exceptions:\newline
\underline{$B_2\rightarrow C_3$}:
\begin{diagram}[size=1em,bottom]
\stackrel{1}{\bullet}&\rImplies&\stackrel{2}{\bullet}&\rLine&\stackrel{\star}{\bullet}
\end{diagram}
In this case $\alpha = \alpha_2$ and $l=\frac{1}{2}$ and thus there is no weight extending $B_2$ to $C_3$.

\underline{$B_3\rightarrow F_4$}:
\begin{diagram}[size=1em,bottom]
\stackrel{1}{\bullet}&\rLine&\stackrel{2}{\bullet}&\rImplies&\stackrel{3}{\bullet}&\rLine&\stackrel{\star}{\bullet}
\end{diagram}
Here  $\alpha = \alpha_3$ and $l = \frac{1}{2}$ and so there is no weight extending $B_3$ to $F_4$.\par
It remains to show that the data from the table leads to Nichols algebras of finite Gelfand-Kirillov dimension. It is clear that the braiding is of exponential type in every case. Moreover in each line the extended Cartan matrix $(b_{ij})$ is of finite type and thus the Nichols algebra has finite Gelfand-Kirillov dimension by theorem \ref{thm_finkrit}.
\qed\end{pf}

\begin{rem}
In the case $A_n\rightarrow A_{n+1}$ the braiding is of Hecke type and this example is already treated in \cite{Rosso_invent}. The cases $A_1\rightarrow A_@$, $A_1\rightarrow B_2$ and $A_1\rightarrow G_2$ (using $\lambda_1$) are treated in \cite{Ich_PBW}, where the relations and the PBW basis are calculated explicitly.
\end{rem}

\section{Results on relations}
\label{sect_rel}
This section contains results on the defining relations of the Nichols algebras studied in section \ref{sect_uqgapp}. First consider a more general setting again.\par
Let $H$ be a Hopf algebra with bijective antipode, $V\in{}^H_H\mYD$, $A:=\Nichols(V)\# H$ and $M$ a graded Yetter-Drinfel'd module over $A$ of highest weight. Furthermore $T_c(M)$ resp. $T_c(M_H\oplus V)$ denote the tensor algebras of $M$ resp. $M_H\oplus V$ viewed as braided Hopf algebras in the corresponing Yetter-Drinfel'd categories ${}^A_A\mYD$ resp. ${}^H_H\mYD$.
The following diagram of $H$-linear maps describes the situation of this section. 

\begin{diagram}
T_c(M)&\pile{\rInto\\\lOnto_{_{T_c(M)\#\eps}}}&T_c(M)\# \Nichols(V)&\lOnto^{_p}&T_c(M_H\oplus V)\\
\dOnto<{_\pi}&&\dOnto<{_{\pi\#\Nichols(V)}}&&\dOnto>{_q}\\
\Nichols(M)&\pile{\rInto\\\lOnto_{_{\Nichols(M)\#\eps}}}&\Nichols(M)\#\Nichols(V)&\rTo^\isom&\Nichols(M_H\oplus V)
\end{diagram}
The maps $p$ and $q$ are the unique algebra morphisms (and braided bialgebra morphisms) that restrict to the identity on $M_H\oplus V$. All maps but $T_c(M)\#\eps$ and $\Nichols(M)\#\eps$ are algebra morphisms. The following proposition is the central tool of this section.

\begin{prop}
\label{prop_idgen}
Assume the situation described above. Fix a subset $X\subset T_c(M_H\oplus V)$ such that $H\cdot X$ generates $\ker q$ as an ideal. Furthermore write the elements of $p(X)$ in the form
\[ p(x) = \sum\limits_i m_i^x\#v_i^x \in T_c(M)\#\Nichols(V)\]
with $m_i^x\in T_c(M),v_i^x\in\Nichols(V)$.
Consider the space
\[ \hat{X} := \left\{ \sum_i m^x_i((v^x_i\#1)\cdot m)|x\in X,m\in T_c(M)\right\}\subset T_c(M).\]
Then $A\cdot \hat{X}$ generates $\ker \pi$ as an ideal.
\end{prop}
\begin{pf}
Let $I:=\ker q$. Obviously
\[ \ker (\pi\#\Nichols(V)) = p(I),\]
and it is easy to check that this implies
\[\ker\pi = (T_c(M)\#\eps)p(I). \]
As $H\cdot X$ generates $I$ as an ideal, $H\cdot p(X)$ generates $p(I)$ as an ideal. Now $T_c(M)\#\eps$ is in general not an algebra morphism, so it is not easy to find ideal generators for $\ker\pi$. The elements of the form
\[ (m'\#v)(h\cdot p(x))(m\#v') = h\s2 \cdot\left[\left(S^{-1}(h\s1)\cdot(m'\#v)\right)\,p(x)\,\left(S(h\s3)\cdot(m\#v')\right)\right]\]
($m,m'\in T_c(M),v,v'\in\Nichols(V),x\in X,h\in H$) generate $p(I)$ as a vector space. Thus $p(I)$ is generated as $H$-module by elements of the form
\begin{eqnarray*}
&&(m'\#v)p(x)(m\#v') = \\
&&=\sum\limits_{i}m'\left((v\su1\#v\su2\sm1)\cdot\left(m_i^x\left((v_i^x{}\su1\#v_i^x{}\su2\sm1)\cdot m\right)\right)\right)\#v\su2\s0v_i^x{}\su2\s0v'
\end{eqnarray*}
($m,m'\in M,v,v'\in\Nichols(V),x\in X$). Now apply the $H$-linear map $T_c(M)\#\eps$ and obtain $H$-module generators of $\ker\pi$ of the form
\[ \sum\limits_{i}m'\left((v\#1)\cdot\left(m_i^x\left((v_i^x\#1)\cdot m\right)\right)\right)\]
($m,m'\in T_c(M),v\in\Nichols(V),x\in X$). Using that $T_c(M)$ is an $H$-module algebra conclude that elements of the form
\[ m'\left((v\#h)\cdot\left(\sum\limits_i m_i^x\left((v_i^x\#1)\cdot m\right)\right)\right)\]
($m,m'\in T_c(M),v\in\Nichols(V),h\in H,x\in X$) generate $\ker\pi$ as vector space. This means that $A\cdot\hat{X}$ generates $\ker\pi$ as (left) ideal in $T_c(M)$.
\qed\end{pf}

\subsection{The quantum group case}
The description of the generators of the ideal $\ker\pi$ obtained in the preceeding theorem is not very explicit as the set $\hat{X}$ may be very large. Nevertheless it is sufficient for the case treated in section \ref{sect_uqgapp}. Assume that $k$ is an algebraically closed field of characteristic zero. Let $M$ be a finite-dimensional integrable $U_q(\mathfrak{g})$-module as in section \ref{sect_uqgmod} with a braiding $c^f$ of strong exponential type with function $\varphi$; moreover assume that the extended Cartan matrix is a generalized Cartan matrix. Let $\hat{U} = \Nichols(V)\#kG$ be the extension defined in \ref{sect_ydstruct}. Furthermore assume that the ideal $\ker q$ is generated by the quantum Serre relations
\begin{eqnarray*}
\forall\alpha,\beta\in\Pi,\alpha\neq\beta: && r_{\alpha\beta} = \ad_c(\hat{F}_\alpha)^{1-b_{\alpha\beta}}(\hat{F}_\beta),\\
\forall \alpha\in\Pi,1\leq i\leq r: &&r_{i\alpha} = \ad_c(m_i)^{1-b_{i\alpha}}(\hat{F}_\alpha),\\
\forall \alpha\in\Pi,1\leq i\leq r: &&r_{\alpha i} = \ad_c(\hat{F}_\alpha)^{1-b_{\alpha i}}(m_i)),\\
\forall 1\leq i\neq j\leq r: &&r_{ij} = \ad_c(m_i)^{1-b_{ij}}(m_j).
\end{eqnarray*}

\begin{rem}
Note that if $\varphi$ is symmetric, then $(b_{ij})$ is a symmetrizeable generalized Cartan matrix and the braiding on $M_{kG}\oplus V$ is of DJ-type by the proof of theorem \ref{thm_finkrit}. In this case \cite[Theorem 2.9]{AS4} ensures that $\ker q$ is generated by these relations.
\end{rem}

In order to apply proposition \ref{prop_idgen} calculate the images of these elements under $p$. First observe
\begin{eqnarray*}
p(r_{\alpha\beta}) &=& \ad_c(p(\hat{F}_\alpha))^{1-b_{\alpha\beta}}(p(\hat{F}_\beta)) = \ad_c(1\#\hat{F}_\alpha)^{1-b_{\alpha\beta}}(1\#\hat{F}_\beta) = \\
 &=&  1\#\left(ad_c(\hat{F}_\alpha)^{1-b_{\alpha\beta}}(\hat{F}_\beta)\right) = 0
\end{eqnarray*}
because this is a relation in $\Nichols(V)$. For $r_{i\alpha}$ use the explicit form of the quantum Serre relations from \cite[Equation A.8]{AS2}:
\[ \ad_c(x)^n(y) = \sum\limits_{s=0}^n (-1)^s\qgauss{n}{s}{\gamma} \gamma^{\frac{s(s-1)}{2}} \eta^s x^{n-s}y x^s,\]
if $c(x\otimes y)= \eta y\otimes x$ and $c(x\otimes x)=\gamma x\otimes x$.
Define the coefficients $q_{xy},x,y\in P$ as in the proof of theorem \ref{thm_finkrit}. Then
\begin{eqnarray*}
p(r_{i\alpha}) &=& p\left(\sum\limits_{s=0}^{1-b_{i\alpha}} (-1)^s\qgauss{1-b_{i\alpha}}{s}{q_{ii}} q_{ii}^{\frac{s(s-1)}{2}} q_{i\alpha}^s m_i^{1-b_{i\alpha}-s}\hat{F}_\alpha m_i^s\right)\\
&=&\sum\limits_{s=0}^{1-b_{i\alpha}} (-1)^s\qgauss{1-b_{i\alpha}}{s}{q_{ii}} q_{ii}^{\frac{s(s-1)}{2}} q_{i\alpha}^s (m_i^{1-b_{i\alpha}-s}\#1)(1\#\hat{F}_\alpha) (m_i^s\#1)\\
&=&\phantom{+}\sum\limits_{s=0}^{1-b_{i\alpha}} (-1)^s\qgauss{1-b_{i\alpha}}{s}{q_{ii}} q_{ii}^{\frac{s(s-1)}{2}} q_{i\alpha}^s \left(m_i^{1-b_{i\alpha}-s}(K_\alpha\cdot m_i^s)\#\hat{F}_\alpha\right)\\
&&+\sum\limits_{s=0}^{1-b_{i\alpha}} (-1)^s\qgauss{1-b_{i\alpha}}{s}{q_{ii}} q_{ii}^{\frac{s(s-1)}{2}} q_{i\alpha}^s \left(m_i^{1-b_{i\alpha}-s}(\hat{F}_\alpha\cdot m_i^s)\#1\right)\\
\end{eqnarray*}
The first summand is zero. This can be seen using
\[ q_{i\alpha}q_{\alpha i} = q_{ii}^{b_{i\alpha}}\:\:\mbox{and}\]
\begin{eqnarray*}
&&\sum\limits_{s=0}^{1-b_{i\alpha}} (-1)^s\qgauss{1-b_{i\alpha}}{s}{q_{ii}} q_{ii}^{\frac{s(s-1)}{2}} q_{i\alpha}^s \left(m_i^{1-b_{i\alpha}-s}(K_\alpha\cdot m_i^s)\#\hat{F}_\alpha\right)\\
&=& \sum\limits_{s=0}^{1-b_{i\alpha}} (-1)^s\qgauss{1-b_{i\alpha}}{s}{q_{ii}} q_{ii}^{\frac{s(s-1)}{2}} q_{i\alpha}^s q_{\alpha i}^s \left(m_i^{1-b_{i\alpha}}\#\hat{F}_\alpha\right)\\
&=&\left(\sum\limits_{s=0}^{1-b_{i\alpha}} (-1)^s\qgauss{1-b_{i\alpha}}{s}{q_{ii}} q_{ii}^{\frac{s(s+1)}{2}} q_{ii}^{-s(1-b_{i\alpha})} \right)\left(m_i^{1-b_{i\alpha}}\#\hat{F}_\alpha\right)=0,\\
\end{eqnarray*}
where the first sum in the last equation is zero by \cite[Equation A.5]{AS2}.
Thus the image of $r_{i\alpha}$ is
\[ p(r_{i\alpha}) = \left(\sum\limits_{s=0}^{1-b_{i\alpha}} (-1)^s\qgauss{1-b_{i\alpha}}{s}{q_{ii}} q_{ii}^{\frac{s(s-1)}{2}} q_{i\alpha}^s m_i^{1-b_{i\alpha}-s}(\hat{F}_\alpha\cdot m_i^s)\right)\#1.\]
$r_{\alpha i}$ is mapped to
\begin{eqnarray*}
p(r_{\alpha i}) &=& \ad_c(p(\hat{F}_\alpha))^{1-b_{\alpha i}}(p(m_i))
= \ad_c(1\#\hat{F}_\alpha)^{1+\frac{2(\alpha,\lambda_i)}{(\alpha,\alpha)}}(m_i\#1)=\\
&=& \left(\hat{F}_\alpha^{1+\frac{2(\alpha,\lambda_i)}{(\alpha,\alpha)}}\cdot m_i\right)\#1 = 0
\end{eqnarray*}
because the braided adjoint action of $\Nichols(V)$ on $M$ (in $T_c(M)\#\Nichols(V)$) is the same as the module action (denoted by $\cdot$) of $\Nichols(V)\subset U_q(\mathfrak{g})$ on $M$. By \cite[5.4.]{Jantzen} the last equality holds. This leaves $r_{ij}$ to be considered. 

\begin{eqnarray*}
p(r_{ij}) &=& p\left(\ad_c(m_i)^{1-b_{ij}}(m_j)\right) 
= \ad_c\left(p(m_i)\right)^{1-b_{ij}}\left(p(m_j)\right)\\
&=& \ad_c(m_i\#1)^{1-b_{ij}}(m_j\#1) = \ad_c(m_i)^{1-b_{ij}}(m_j)\#1
\end{eqnarray*}

\begin{rem}
A short calculation results in the following representation of the relation coming from $r_{i\alpha}$ for $1\leq i\leq r,\alpha\in\Pi$:
\[ R_{i\alpha} := \sum\limits_{t=0}^{1-b_{i\alpha}}\left[\sum\limits_{s=t+1}^{1-b_{i\alpha}}(-1)^s \qgauss{1-b_{i\alpha}}{s}{q_{ii}} q_{ii}^{\frac{s(s-1)}{2}}q_{ii}^{sb_{i\alpha}}\right] q_{\alpha i}^{-1-t} m_i^{-b_{i\alpha}-t}\left(\hat{F}_\alpha \cdot m_i\right) m_i^t. \]
The relations coming from the $r_{ij},1\leq i,j\leq r$ are
\[
R_{ij} := \ad_c(m_i)^{1-b_{ij}}(m_j)
\]
It follows from proposition \ref{prop_idgen} that the $\Nichols(V)\#kG$ submodule of $T_c(M)$ generated by the elements
\[ \{R_{i\alpha}|1\leq i\leq r,\alpha\in\Pi\}\cup\{R_{ij}|1\leq i\neq j\leq r\}\]
generates the kernel of the canonical map
\[ \pi:T_c(M)\rightarrow \Nichols(M)\]
as an ideal.
\end{rem}

\begin{thm}
Let $M$ be a finite-dimensional integrable $U_q(\mathfrak{g})$-module and fix a braiding $c^f$ as in section \ref{sect_uqgmod} of strong exponential type with symmetric function $\varphi$. Assume that the extended Cartan matrix $(b_{ij})_{i,j\in P}$ is a generalized Cartan matrix. Consider the grading on $\Nichols(M)$ such that the elements of $M$ have degree 1.
Then $\Nichols(M)$ is generated by $M$ with homogenous relations of the degrees
\[ 2-b_{ij}\:\mtxt{for}\: 1\leq i\neq j\leq r \:\:\mbox{and}\]
\[ 1-b_{i\alpha}\:\mtxt{for}\: 1\leq i\leq r,\alpha\in\Pi\:\:\mbox{such that}\:\:b_{i\alpha}\neq 0.\]
\end{thm}
The last column of table \ref{table_extensions} was calculated using this theorem.
\begin{pf}
Again realize the module $M$ as a Yetter-Drinfel'd module over $\Nichols(V)\#kG$ as in section \ref{sect_ydstruct}. The $\Nichols(V)\#kG$-module generated by the elements $R_{i\alpha},R_{ij}, 1\leq i\neq j\leq r,\alpha\in\Pi$ generates the $\ker\pi$ as an ideal. $R_{i\alpha}$ has degree $1-b_{i\alpha}$, $R_{ij}$ has degree $2-b_{ij}$ (with respect to the grading of $T_c(M)$ giving $M$ the degree 1).
As the homogenous components of $T_c(M)$ are $\Nichols(V)\#kG$-modules all defining relations can be found in the degrees
\[ 1-b_{i\alpha}\:\:\mbox{and}\:\:2-b_{ij}.\]
Observe that $p(r_{i\alpha})$ is zero if $b_{i\alpha}=0$: The summand for $s=0$ is zero anyway because $\eps(\hat{F}_\alpha)=0$. The summand for $s=1$ is a scalar multiple of $\hat{F}_\alpha\cdot m_i$. If $b_{i\alpha}=0$ then also $b_{\alpha i}=0$ and thus by \cite[5.4.]{Jantzen} $\hat{F}_\alpha\cdot m_i = 0$.
\qed\end{pf}

\bibliography{../../tex/promotion.bib}
\end{document}